\documentclass[a4paper,10pt]{article}

\usepackage{amsmath,amsfonts,amssymb,amsthm,epsfig}

\def\C{\mathbb{C}}

\def\N{\mathbb{N}}

\def\sl{\mathfrak{sl}_2}

\DeclareMathOperator{\im}{Im} 
\DeclareMathOperator{\Hom}{Hom}
\DeclareMathOperator{\rank}{rank}
\DeclareMathOperator{\End}{End}

\newtheorem*{theo}{Theorem}
\newtheorem*{prop}{Proposition}

\numberwithin{equation}{section}

\newcommand{\id}{\ensuremath{\mathbf {1}}}

\newcommand{\Span}{\operatorname{Span}}

\begin{document}
\title{Quiver varieties and fusion products for $\mathfrak{sl}_2$}
\author{Alistair Savage and Olivier Schiffmann}
\date{}
\maketitle
\noindent

\paragraph{Introduction.} In a remarkable series of work starting in \cite{N0},
Nakajima gives a geometric realization of integrable highest weight 
representations $V_\lambda$ of a Kac-Moody algebra $\mathfrak{g}$ in the
homology of a certain Lagrangian subvariety $\mathcal{L}(\lambda)$ of a
symplectic variety $\mathcal{M}(\lambda)$ constructed from the Dynkin
diagram of $\mathfrak{g}$ (the \textit{quiver variety}). In particular, in
\cite{N2}, he realizes the tensor product $V_\lambda \otimes V_{\mu}$
as the homology of a ``tensor product variety'' $\mathcal{L}(\lambda,
\mu) \subset \mathcal{M}(\lambda + \mu)$ (the same construction also appears 
independently in \cite{M}). When
$\mathfrak{g}$ is simple, one might ask if a similar construction can produce
the \textit{fusion} tensor products $V_{\lambda} \otimes_l V_{\mu}$,
certain truncations of $V_\lambda \otimes V_\mu$.

\paragraph{}In this short note, we answer this question affirmatively when
$\mathfrak{g}=\mathfrak{sl}_2$. In this case, $V_{\lambda}
\otimes_l V_{\mu}$ is realized as the homology of the most natural
subvarieties $\mathcal{L}_l(\lambda,\mu) \subset \mathcal{L}(\lambda,\mu)$
(see Section 3). We also consider the case of a tensor product of
arbitrarily many $\mathfrak{sl}_2$-modules $V_{\lambda_1}, \cdots ,V_{\lambda_r}$. Finally, we
give a combinatorial description of the irreducible components
of $\mathcal{L}_l(\lambda,\mu)$ (and $\mathcal{L}_l(\lambda_1, \ldots, 
\lambda_r)$) using the notions of graphical
calculus and crossingless matches for $\mathfrak{sl}_2$ (see \cite{FK97} 
and \cite{Savage}). We do not expect these constructions to generalize to Lie
algebras of higher rank.

\centerline{\textbf{Acknowledgments.}}
We would like to thank I. Frenkel for raising the question of the relation
between quiver varieties and fusion products. The research of the first
author was supported in part by the Natural Science and Engineering
Research Council of Canada. 


\section{Fusion products for $U(\mathfrak{sl}_2)$.}
\paragraph{1.1.} Let 
$\mathcal{R}$ denote the category of finite-dimensional $\sl$-modules, and 
for $i \geq 0$ let $V_i$ denote the simple module of highest weight $i$. Let 
$\C[\mathcal{R}]$ be the Grothendieck ring of $\mathcal{R}$ and let $[V]$ 
denote the class of a module $V$. We have
$$V_i \otimes V_j \simeq \bigoplus_{k=j-i}^{i+j} V_k,\qquad [V_i] \cdot [V_j]
=\sum_{k=j-i}^{i+j} [V_k],\qquad \mathrm{for\;} i \leq j$$
where in the sums $k$ increases by twos.
\paragraph{1.2.} Now let us fix some positive integer $l \in \N$. Consider the
quotient
$$\C_l[\mathcal{R}]=\C[\mathcal{R}]/[V_{l+1}]\C[\mathcal{R}].$$
Denoting by $[V]_l$ the image of $[V]$ in $\C_l[\mathcal{R}]$, we have
 $\C_l[\mathcal{R}]=\C[V_0]_l \oplus \cdots \oplus \C[V_l]_l$, and
$$[V_i \otimes V_j]_l=\sum_{k=j-i}^{\min(i+j,2l -i-j)} [V_k]_l,\qquad 
\mathrm{for\;} 
i \leq j \leq l.$$
We also set
$$V_i \otimes_l V_j=\bigoplus_{k=j-i}^{\min(i+j,2l-i-j)} V_k,\qquad
\mathrm{for\;} 
i \leq j \leq l.$$
Again, in the above sums, $k$ increases by twos.
The ring $\C_l[\mathcal{R}]$ appears in conformal field theory (as the 
Grothendieck ring of the modular category of integrable 
$\widehat{\mathfrak{sl}}_2$-modules of level $l$) and in quantum group theory
(as the Grothendieck ring of a suitable quotient of the category of tilting
modules over $U_\epsilon(\sl)$ when $\epsilon$ is an $l$th root of unity).


\section{Lagrangian construction of $U(\sl)$.}
We briefly recall Ginzburg's construction of irreducible representations of
$\sl$ in the homology of certain varieties associated to partial flag varieties
(cf. \cite{Ginz}). 
We use the (in this case equivalent) language of quiver varieties 
(cf. \cite{Nak}).
\paragraph{2.1.} Let $v,w \in \N$ and let $V$ and $W$ be $\C$-vector spaces
of dimensions $v$ and $w$ respectively. Consider the space
$$M(v,w)=\{(i,j)\;|\;ij=0; \;\ker j=\{0\}\} \subset 
\mathrm{Hom}\;(W,V) \oplus \mathrm{Hom}\;(V,W).$$
We let $GL(V)$ act on $M(v,w)$ via $g\cdot(i,j)=(gi,jg^{-1})$. This action is
free and we set $\mathcal{M}(v,w)=M(v,w)/GL(V)$. The assignment
$(i,j) \mapsto (ji, \mathrm{Im}\;j)$ defines an isomorphism between $\mathcal{M}(v,w)$ and the variety
$$\mathcal{F}_{v,w}=\{(t,V_0)\;| V_0 \subset W, \mathrm{dim}\;V_0=v,
\mathrm{Im}\;t \subset V_0 \subset \ker t \} \subset
\mathcal{N}_W \times \mathrm{Gr}(v,w),$$
where $\mathcal{N}_W$ is the nullcone of $\mathfrak{gl} (W)$ and 
$\mathrm{Gr}(v,w)$ is the Grassmannian of $v$-dimensional subspaces in $W$.
We will denote by $\pi: \mathcal{M}(v,w) \to \mathcal{N}_W$, the
projection $(i,j) \mapsto ji$. For any $t \in \mathcal{N}_W$ such that $t^2=0$ we set
$\mathcal{M}(v,w)_t=\pi^{-1}(t)$ and $\mathcal{M}(w)_t = \sqcup_v
\mathcal{M}(v,w)_t$. In particular, we set $\mathcal{L}(v,w)=
\pi^{-1}(0)$. Observe that $\mathcal{L}(v,w)$ is just $\mathrm{Gr}(v,w)$ and
that $\mathcal{M}(v,w)$ is isomorphic to the cotangent bundle of
$\mathcal{L}(v,w)$.
We have $\mathrm{dim}\;\mathcal{M}(v,w)=2 
\mathrm{dim}\;\mathcal{L}(v,w)= 2v(w-v)$. 
For $v_1,v_2,w \in \N$ we also consider the variety of 
triples
$$Z(v_1,v_2,w)=\{((i_1,j_1),(i_2,j_2))\;|\;j_1i_1=j_2i_2\}\subset 
\mathcal{M}(v_1,w) \times \mathcal{M}(v_2,w). $$
Then $\mathrm{dim}\;Z(v_1,v_2,w)=v_1(w-v_1)+v_2(w-v_2)$.
\paragraph{}The form $\omega((i,j),(i',j'))=\mathrm{Tr}_V(ij'-i'j)$ defines a
symplectic structure on $\mathcal{M}(v,w)$, for which the variety 
$\mathcal{L}(v,w)$ is Lagrangian. Equip $\mathcal{M}(v_1,w) \times 
\mathcal{M}(v_2,w)$ with the symplectic form $\omega \times (-\omega)$. Then 
$Z(v_1,v_2,w)$ is also Lagrangian.  Let $Z(w) = \sqcup_{v_1,v_2} Z(v_1,v_2,w)$.

\paragraph{2.2.} For any complex algebraic variety $X$ we let $H_*(X)$ be the
Borel-Moore homology with coefficients in $\C$, and set $H_{top}(X)=H_{2d}(X)$
where $d=\mathrm{dim}\;X$.
\paragraph{} Let $p_{ij}: \mathcal{M}(v_1,w) \times \mathcal{M}(v_2,w)
\times \mathcal{M}(v_3,w) \to \mathcal{M}(v_i,w) \times \mathcal{M}(v_j,w)$
be the obvious projections. The map
$$p_{13}:\;p_{12}^{-1}(Z(v_1,v_2,w)) \cap p_{23}^{-1}(Z(v_2,v_3,w)) \to
Z(v_1,v_3,w)$$
is proper and we can define the convolution product
\begin{align*}
H_i(Z(v_1,v_2,w)) \otimes H_j(Z(v_2,v_3,w)) &\to H_{i+j-d_2}(Z(v_1,v_3,w))\\
c \otimes c' &\mapsto p_{13*}(p_{12}^*(c) \cap p_{23}^*(c'))
\end{align*}
where $d_2=4v_2(w-v_2)$. In particular, this gives rise to an algebra 
structure on $H_{top}(Z(w))=
\bigoplus_{v_1,v_2} H_{top}(Z(v_1,v_2,w))$.

\paragraph{}Now let $t \in \mathcal{N}_W$ such that $t^2=0$. The projection
$$p_1:\;Z(v_1,v_2,w) \cap p_2^{-1}(\mathcal{M}(v_2,w)_t) \to 
\mathcal{M}(v_1,w)_t$$
(where $p_1$ and $p_2$ are the obvious projections)
is proper and the convolution action
\begin{align*}
H_{top}(Z(v_1,v_2,w)) \otimes H_{top}(\mathcal{M}(v_2,w)_t) &\to H_{top}
(\mathcal{M}(v_1,w)_t)\\
c \otimes c' &\mapsto p_{1*}(c \cap p_{2}^*(c'))
\end{align*}
makes $H_{top}(\mathcal{M}(w)_t)=\bigoplus_v H_{top}(\mathcal{M}(v,w)_t)$
into a $H_{top}(Z(w))$-module.

\begin{theo}[\cite{Ginz}] There is a natural surjective homomorphism
$\Phi:\;U(\sl) \to H_{top}(Z(w))$. Under $\Phi$, the module 
$H_{top}(\mathcal{M}(w)_t)$ is isomorphic to $V_{w-2u}$ where 
$u=\rank t$.\end{theo}

\paragraph{2.3.} We now give the realization of tensor products of 
$U(\sl)$-modules. Let $w=w_1+\cdots + w_r$ and fix $W=W_1 \oplus \cdots \oplus
W_r$ with $\mathrm{dim}\;W_i=w_i$. Let $W_0 = 0$. The group $GL(W)$ acts on 
$\mathcal{M}(v,w)$ by $g\cdot (i,j)=(ig^{-1},gj)$. Consider the embedding
\begin{align*}
\sigma:\;(\C^*)^{r-1} &\to \prod_{i=1}^r GL(W_i) \subset GL(W)\\
(t_2,t_3,\ldots,t_r) &\mapsto (Id, t_2^{-1}, t_2^{-1}t_3^{-1},\ldots,
t_2^{-1} \cdots t_r^{-1})
\end{align*}
Then, for each $v$, we have (see e.g \cite[Lemma 3.2]{N2})
$$\mathcal{M}(v,w)^\sigma =\underset{v_1+\cdots + v_r=v}{\bigsqcup}
\mathcal{M}(v_1,w_1) 
\times \cdots \times \mathcal{M}(v_r,w_r).$$
Consider the subvarieties
\[
\mathcal{M}(v,w_1,\ldots,w_r)=\{x \in \mathcal{M}(v,w)\;|\;
\underset{t_i \to 0}{\mathrm{lim}} \sigma(t_2,\ldots,t_r)\cdot x
\;\mathrm{exists}\}
\]
$$\mathcal{N}_W(w_1,\ldots,w_r)=\{t \in \mathcal{N}_W\;|\;
\underset{t_i \to 0}{\mathrm{lim}} \sigma(t_2,\ldots,t_r)\cdot t
\;\mathrm{exists}\}.
$$
For $x \in \mathcal{M}(v,w_1,\ldots,w_r)$, let us set $\tau(x)=
\underset{t_i \to 0}{\mathrm{lim}} \sigma(t_2,\ldots,t_r)\cdot x$.
We define $\tau(t)$ similarly for $t \in \mathcal{N}_W(w_1,\ldots,w_r)$.  
Now 
consider
$$\mathcal{L}(v,w_1,\ldots,w_r)=\{x \in \mathcal{M}(v,w_1,\ldots,w_r)\;|\;
\tau(x)
\in \prod_i \mathcal{L}(v_i,w_i)\;\mathrm{for\;some\;}(v_i)\}.
$$
Set $\mathcal{L}(w_1,\ldots,w_r)=\sqcup_v \mathcal{L}(v,w_1,\ldots,w_r)$.
Note that $\mathcal{L}(w_1,\ldots,w_r)=\pi^{-1} (\tau^{-1}(0))$
so that we have an action of $H_{top}(Z(w))$ on 
$H_{top}(\mathcal{L}(w_1,\ldots,w_r))$.
Moreover, it is easy to check that 
$\mathcal{L}(w_1,\ldots,w_r)$ is Lagrangian.
Note that $\mathcal{L}(w_1,\ldots,w_r)$ is isomorphic to the variety
\[
\{(t,V_0) \;|\; V_0 \subset W,\, \im t \subset V_0 \subset \ker t,\,
t(W_j) \subset W_0 \oplus \dots \oplus W_{j-1},\, 1\le j \le r\}.
\]

\begin{theo}[\cite{GRV}, \cite{N2}, \cite{M}]  
$H_{top}(\mathcal{L}(w_1,\ldots,w_r))$ is isomorphic
to $V_{w_1} \otimes \cdots \otimes V_{w_r}$ as a $U(\sl)$-module.
\end{theo}


\section{Lagrangian construction of the fusion product}
\paragraph{}Let us fix some positive integer $l$. We will now describe 
an open subvariety of $\mathcal{L}(w_1,\ldots,w_r)$ whose homology realizes
the fusion product $V_{w_1} \otimes_l \cdots \otimes_l V_{w_r}$.

\paragraph{3.1.} We keep the notation of 2.3.
For all $k \in \N$ and $t \in \mathcal{N}_{W_1 \oplus \cdots \oplus
W_k} (w_1,\ldots, w_k)$ we set 
$\tau_k(t)=\mathrm{lim}_{t_k\to 0} \sigma (1,\ldots,1,t_k) (t)$.
Let us consider the open subvariety $\mathcal{N}^l(w_1,w_2)=
\{t \in \mathcal{N}_{W_1 \oplus W_2}\;|
\dim \ker t \leq l \}$ of $\mathcal{N}_{W_1 \oplus W_2}$ and define
inductively
\begin{equation}
\label{Nl_def}
\begin{split}
\mathcal{N}^l(w_1,\ldots,w_k)=\{t \in \mathcal{N}_{W_1 \oplus \cdots \oplus
W_k}\;|\dim \ker t \leq l + \rank \tau_k(t), \\
t|_{W_1 \oplus \dots \oplus W_{k-1}} \in \mathcal{N}^l(w_1,\dots,w_{k-1}) \}
\end{split}
\end{equation}
for $k \ge 3$.
Finally, set $\mathcal{L}_l(w_1,\ldots,w_r)=\mathcal{L}(w_1,\ldots,w_r)
\cap \pi^{-1} (\mathcal{N}^l(w_1,\ldots,w_r))$.
By definition $\mathcal{L}_l(w_1,\ldots,w_r)$ is an open subvariety of 
$\mathcal{L}(w_1,\ldots,w_r)$ and therefore
$H_{top}(\mathcal{L}_l(w_1,\ldots,w_r))$
is a $H_{top}(Z(w))$-module.

\begin{theo} $H_{top}(\mathcal{L}_l(w_1,\ldots,w_r))$ is isomorphic to
$V_{w_1} \otimes_l \cdots \otimes_l V_{w_r}$ as a $U(\sl)$-module.
\end{theo}
\noindent
\begin{proof} We proceed by induction. Suppose $r=2$. It is enough to 
describe the irreducible components of $\mathcal{L}_l(w_1,w_2)$ corresponding
to highest weight vectors in the $U(\sl)$-module 
$H_{top}(\mathcal{L}_l(w_1,w_2))$. The irreducible components of 
$\mathcal{L}(w_1,w_2)$ corresponding to highest-weight vectors are
the
$$I_v=\{(i,j)\;|\;j(V) \subset W_1,\;i(W_2)=V,\;i(W_1)=0\},
\qquad\mathrm{for\;} 0 \leq v \leq w_1,w_2$$
and the associated highest weight is $w_1+w_2-2v$. Note that the condition
$\dim \ker ji \leq l$ is equivalent to the condition $w_1+ w_2-2v 
\leq 2l - w_1-w_2$. Now suppose that the theorem is proved for tensor
products of $r-1$ modules, and let us set $W'=W_1 \oplus\cdots\oplus W_{r-1}$.
For each $u \in \N$ let us set $\mathcal{N}_{W'}(u)=\{t \in \mathcal{N}_{W'}\;
| \rank t=u\}$. Recall that $\mathcal{L}_l(w_1,\ldots,w_{r-1})$
is Lagrangian and that $\pi$ is semi-small with all strata being relevant
(c.f \cite[$\S10$]{Nak}). Thus $\pi(\mathcal{L}_l(w_1,\ldots ,w_{r-1})) \cap
\mathcal{N}_{W'}(u)$ is a subvariety of $\mathcal{N}_{W'}(u)$ of dimension
$\frac{1}{2}\mathrm{dim\;}\mathcal{N}_{W'}(u)$. Let $\mathcal{C}_1^u,\ldots ,
\mathcal{C}_{s(u)}^u$ be its irreducible components. By the induction 
hypothesis,
\begin{equation}\label{E:1}
s(u)=\mathrm{dim\;Hom}_{\sl}(V_{w'-2u},V_{w_1} \otimes_l \cdots \otimes_l
V_{w_{r-1}}).
\end{equation}
The irreducible components of $\mathcal{L}_l(v,w_1,\ldots,w_r)$ 
corresponding to highest weight vectors of $H_{top}(\mathcal{L}_l(w_1,\ldots,
w_r))$ are of the form $\overline{I_\chi}$ with
$$I_{\chi}=\{(i,j)\;|i(W)=V,\;j(V) \subset W',\;(i_{W'},j) \in \chi\}$$
where $\chi$ is an irreducible component of $\mathcal{L}_l
(v,w_1,\ldots,w_{r-1})$, and the associated highest weight is $w-2v$ (note that
$I_\chi$ may be empty). Let us fix $u \in \N$ and $\mathcal{C}_k^u$ for some 
$k \leq s(u)$. Let 
$\chi \subset \overline{\pi^{-1}(\mathcal{C}_k^u)}
 \cap \mathcal{L}_l(v,w_1,\ldots,w_{r-1})$
be an irreducible component. Then $I_\chi \subset \overline{\mathcal{L}_l(w_1,\ldots,
w_r)}$ if for all $(i,j)$ in (an open dense subset of) $I_{\chi}$ we have
$\mathrm{dim\;Im\;}ji \leq l +u.$ This is equivalent to the condition that
the corresponding highest weight $w-2v$ satisfies
\begin{equation}\label{E:2}
w-2v \leq 2l-w_r-(w'-2u).
\end{equation}
Equations (\ref{E:1}) and (\ref{E:2}) together imply that 
$$H_{top}(\mathcal{L}_l(w_1,\ldots,w_r))\simeq 
(V_{w_1} \otimes_l \cdots 
\otimes_l V_{w_{r-1}})\otimes_l V_{w_r}$$
 as a $U(\sl)$-module, as desired.
\end{proof}

\paragraph{Remarks.} i) The above construction is not canonical in the sense
that it was made using a choice of a bracketing of the tensor product, namely
$$(\cdots (( V_{w_1} \otimes_l V_{w_2}) \otimes_l V_{w_3} )\cdots 
 \otimes_l V_{w_r}).$$
 Different
bracketings give rise to different (possibly non-isomorphic) open
subvarieties of $\mathcal{L}_l(w_1,\ldots,w_r)$ realizing the same
fusion tensor product.\\
ii) One might be tempted to define in an analogous fashion a truncated
tensor product for finite-dimensional representations of 
$U_q(\widehat{\mathfrak{sl}}_2)$ by considering equivariant K-theory
of $\mathcal{L}_l(w_1,w_2)$ rather than Borel-Moore homology. However, it
is easy to check that (because of Remark i)) the resulting product is not 
associative.


\section{A graphical calculus for the fusion product}

\paragraph{4.1.}
We first recall some results on the graphical calculus of tensor
products and intertwiners.  For a more complete treatment, see
\cite{FK97} and \cite{Savage}.
In the graphical calculus, $V_d$ is depicted by a box
marked $d$ with $d$ vertices.  To depict the set $CM_{w_1, \dots,
  w_r}^\mu$ of crossingless matches, we place the
boxes representing the
$V_{w_i}$ on a horizontal line and the box representing
$V_\mu$ on another horizontal line lying above the first one.
$CM_{w_1, \dots, w_r}^\mu$ is then the set of non-intersecting curves
(up to isotopy) connecting the vertices of the boxes such that the following
conditions are satisfied:

\begin{enumerate}
\item Each curve connects exactly two vertices.
\item Each vertex is the end point of exactly one curve.
\item No curve joins a box to itself.
\item The curves lie inside the box bounded by the two horizontal
  lines and the vertical lines through the extreme right and left points.
\end{enumerate}

We call the curves joining two lower boxes \emph{lower curves} and
those joining a lower
and an upper box \emph{middle curves}.  We define the set of oriented
crossingless matches $OCM_{w_1,\dots, w_r}^\mu$ to be the set of
elements of $CM_{w_1, \dots, w_r}^\mu$ along with an
orientation of the
curves such that all lower curves are oriented to the left
and all middle curves are oriented
so that those oriented down are to the right of those oriented up.

As shown in \cite{FK97}, the set of crossingless
matches $CM_{w_1, \dots, w_r}^\mu$ is in
one-to-one correspondence with a basis of the set of intertwiners
\[
H_{w_1, \dots, w_r}^\mu
  \stackrel{\text{def}}{=} \Hom (V_{w_1} \otimes \dots \otimes
  V_{w_r}, V_\mu).
\]
The matrix coefficients of the intertwiner associated to a particular
crossingless match are given by Theorem~2.1 of
\cite{FK97}.

We will also need to define the set of \emph{lower crossingless matches}
$LCM_{w_1, \dots, w_r}^\mu$ and \emph{oriented lower crossingless
  matches} $OLCM_{w_1, \dots, w_r}^\mu$.
Elements of $LCM_{w_1, \dots, w_r}^\mu$ and $OLCM_{w_1, \dots,
  w_r}^\mu$ are obtained from elements of $CM_{w_1, \dots, w_r}^\mu$ and
$OCM_{w_1, \dots, w_r}^\mu$ (respectively) by removing the upper box
(thus converting
lower end points of middle curves to unmatched vertices).
For the case of $OLCM_{w_1, \dots, w_r}^\mu$, unmatched vertices will
still have an orientation (indicated by an arrow attached to the
vertex).  As for middle curves in the case of $OCM_{w_1,\dots, w_r}^\mu$,
the unmatched
vertices in an element of $OLCM_{w_1, \dots, w_r}^\mu$ must be
arranged so that those
oriented down are to the right of those oriented up.

Note that the set of lower crossingless matches
$LCM=LCM_{w_1, \dots, w_r}$ is in one-to-one correspondence with the
set $\bigcup_{\mu} CM_{w_1, \dots, w_r}^\mu$.  From now on, we will
identify these two sets.  

\paragraph{4.2.} Let $s$ be a bracketing of the tensor product $V_{w_1}
\otimes \dots \otimes V_{w_r}$.  Pick an ordering of the tensor
operations compatible with this bracketing.  For
each $n$ such that $1 \le n \le r-1$, let $S_n$ be the set of the
$V_{w_i}$ separated from the $n^{th}$ tensor product operation
only by operations ranked lower than or equal to $n$. Then let
$ ^l_sCM_{w_1, \dots, w_r}^\mu$ be the set of elements of
$CM_{w_1, \dots, w_r}^\mu$ satisfying the following condition:  for
each $n$, the number of curves connecting $V_{w_i}$'s in $S_n$
to either $V_{w_i}$'s in $S_n$ on the other side of the
$n^{th}$ tensor product symbol or $V_{w}$'s not in $S_n$ is
less than or equal to $l$.  Note that this condition does not depend
on the particular ordering so long as it is compatible with the
bracketing $s$.

Let $ ^l_sLCM = {^l_sLCM}_{w_1, \dots, w_r}$ be the set of lower
crossingless matches satisfying the same condition (where unmatched
vertices are always counted as curves with the other end point
outside of any $S_n$) and identify this set with
the set $\bigcup_{\mu} {^l_sCM}_{w_1, \dots,
w_r}^\mu$.  We define $ ^l_sOCM_{w_1, \dots, w_r}^\mu$ and $ ^l_sOLCM =
{^l_sOLCM}_{w_1, \dots, w_r}$ similarly (and the
corresponding identification is made).

Note that in the case $r=2$ the condition in the definition simplifies to the
requirement that the total number of curves (including middle curves)
is less than or equal to $l$.  In fact, the given definition simply
arises from applying this condition to each tensor product operation
(in the given ordering),
neglecting curves with both end points in $V_{w_i}$'s which have
already been tensored together.

\begin{prop}
The set $ ^l_sCM_{w_1, \dots, w_r}^\mu$ is in one-to-one correspondence
with a basis of the space of intertwiners $ ^lH_{w_1, \dots, w_r}^\mu
\stackrel{\text{def}}{=} \Hom (V_{w_1}
  \otimes_l \dots \otimes_l V_{w_r}, V_\mu)$.
\end{prop}

\begin{proof}
We first consider the case $r=2$.  For any $b \in CM_{w_1,
  w_2}^\mu$, the total number of curves is equal to $(w_1
+ w_2  + \mu)/2$ (since each vertex is an end point of exactly one curve).
  Thus the condition that the total number of curves is less than or
  equal to $l$ reduces to $w_1 + w_2 + \mu \le 2l$ or $\mu \le 2l -
  w_1 - w_2$ as desired.

Now assume the result holds for the product of less than $r$
irreducible modules and that for the product of $V_{w_1}$ through
$V_{w_r}$, the $r^{th}$ tensor product operation is the
one occurring between $V_{w_k}$ and $V_{w_{k+1}}$ ($k<r$).  Note that
\[
\bigoplus_\nu {^lH}_{w_1, \dots, w_{k}}^\nu \otimes {^lH}_{\nu,
  w_{k+1}, \dots, w_{r}}^\mu
\cong {^lH}_{w_1, \dots, w_r}^\mu
\]
via the map $f \otimes g \mapsto g(f \otimes
\text{id}_{V_{w_{k+1}} \oplus \dots \oplus V_{w_r}})$.
Now, if $s_1$ is the bracketing of the first $k$ modules and $s_2$ is
the bracketing of the last $r-k$ modules, it is easy to see that
\[
\sum_\nu {^l_{s_1}CM}_{w_1, \dots, w_k}^\nu \times {^l_{s_2}CM}_{\nu, w_{k+1},
  \dots, w_r}^\mu
\cong {^l_sCM}_{w_1, \dots, w_r}^\mu \; \text{(as sets)}.
\]
The result now follows by induction.
\end{proof}

\paragraph{4.3.} From the associativity of the fusion tensor product
it follows immediately that the order of the
set ${^l_sCM}_{w_1, \dots, w_r}^\mu$ is independent
of the bracketing~$s$.  However, we will present
here a direct proof.

\begin{prop}
The order of the set ${^l_sCM}_{w_1, \dots, w_r}^\mu$ is independent
of the bracketing~$s$.
\end{prop}

\begin{proof}
It suffices to prove the statement for three factors.  Let
$s_1$ be the bracketing $(V_{w_1} \otimes V_{w_2}) \otimes
V_{w_3}$ and $s_2$ be the bracketing $V_{w_1} \otimes
(V_{w_2} \otimes V_{w_3})$.
We will set up a one-to-one correspondence between ${^l_{s_1}CM}_{w_1,
\dots, w_r}^\mu$ and ${^l_{s_2}CM}_{w_1, \dots, w_r}^\mu$.  We will
first establish a one-to-one correspondence between the subsets
consisting of those crossingless matches with no curves connecting
$V_{w_1}$ and $V_{w_3}$ and a fixed number $n$ of lower
curves.  Let $a$ (resp. $b$) denote the number of curves connecting
$V_{w_1}$ (resp. $V_{w_3}$) to $V_{w_2}$.  Thus
$a+b = n$.  Now, the number of curves with at least one end point in
$V_{w_1}$ or $V_{w_2}$ is $w_1 + w_2 - a$ and
the total number of curves minus the curves connecting $V_{w_1}$
to $V_{w_2}$ is $w_1 + w_2 + w_3 - n - a$.
Thus a crossingless match lies in ${^l_{s_1}CM}_{w_1,
\dots, w_r}^\mu$ if and only if
\[
w_1 + w_2 - a \le l,\ w_1 + w_2 +
w_3 - n - a \le l.
\]
Similarly, a crossingless match lies in ${^l_{s_2}CM}_{w_1,
\dots, w_r}^\mu$ if and only if
\[
w_2 + w_3 - b \le l,\ w_1 + w_2 +
w_3 - n - b \le l.
\]
Now, the largest possible value of $a$ is $\min(w_1,n)$ and the
largest possible value of $b$ is $\min(w_3,n)$.  Therefore, by
counting the possible values of $a$, the
number of crossingless matches in ${^l_{s_1}CM}_{w_1,
\dots, w_r}^\mu$ with no curves connecting $V_{w_1}$ and
$V_{w_3}$ and with $n$ total curves is equal to
\[
r_a = \min(w_1,n) - \max(w_1 + w_2 - l, w_1
+ w_2 + w_3 - n - l) + 1
\]
if this number is positive and zero otherwise.  Similarly, the
number of crossingless matches in ${^l_{s_2}CM}_{w_1,
\dots, w_r}^\mu$ with no curves connecting $V_{w_1}$ and
$V_{w_3}$ and with $n$ total curves is equal to
\[
r_b = \min(w_3,n) - \max(w_2 + w_3 - l, w_1
+ w_2 + w_3 - n - l) + 1
\]
if this number is positive and zero otherwise.  Considering the four
cases $n \le w_1, w_3$; $n \ge w_1, w_3$;
$w_1 \le n \le w_3$ and $w_3 \le n \le w_1$ we
easily see that $r_a = r_b$ in all cases.

It remains to establish a one-to-one correspondence between the
elements of ${^l_{s_1}CM}_{w_1, \dots, w_r}^\mu$ and
${^l_{s_2}CM}_{w_1, \dots, w_r}^\mu$ with $c \ge 1$ curves joining
$V_{w_1}$ and $V_{w_3}$.  Fix the number of lower curves
with one end point in $V_{w_2}$ to be $n$.  Since $V_{w_1}$
and $V_{w_3}$ are connected, there can be no
middles curves with end points in $V_{w_2}$.  Thus $s=w_2$.
Define $a$ and $b$ as above.  By an argument analogous to that given
in the earlier case, the
number of crossingless matches in ${^l_{s_1}CM}_{w_1,
\dots, w_r}^\mu$ with $c \ge 1$ curves connecting $V_{w_1}$ to
$V_{w_3}$ and with $n$ lower curves with one end point in
$V_{w_2}$ is equal to
\[
r_a = \min(w_1-c,w_2) - \max(w_1 + w_2 - l, w_1
+ w_3 - l - c) + 1
\]
if this number is positive and zero otherwise.  Similarly, the
number of crossingless matches in ${^l_{s_2}CM}_{w_1,
\dots, w_r}^\mu$ with $c \ge 1$ curves connecting $V_{w_1}$ to
$V_{w_3}$ and with $n$ lower curves with one end point in
$V_{w_2}$ is equal to
\[
r_b = \min(w_3-c,w_2) - \max(w_2 + w_3 - l, w_1
+ w_3 - l - c) + 1
\]
if this number is positive and zero otherwise.  Considering the four
cases $w_2 \le w_1-c, w_3-c$; $w_2 \ge
w_1-c, w_3-c$;
$w_1-c \le w_2 \le w_3-c$ and $w_3-c \le
w_2 \le w_1-c$ we
easily see that $r_a = r_b$ in all cases.  This concludes the proof.
\end{proof}

From now on, we will use the bracketing $( \cdots ((V_{w_1}
  \otimes V_{w_2}) \otimes V_{w_3}) \cdots V_{w_r})$
  unless explicitly stated otherwise.  Thus, if we omit a
  subscript $s$, we take $s$ to be this bracketing.


\section{The fusion product via constructible functions}
\paragraph{5.1.}
Fix a $w=w_1 + \dots + w_r$ dimensional $\C$-vector space $W$ and let
\begin{equation*}
\begin{split}
\mathfrak{T}(w_1, \dots, w_r) = \{(\mathbf{D}=\{D_i\}_{i=0}^r,V_0,t)
\;|\; 0 = D_0 \subset D_1 \subset \dots \subset D_r = W,\, V_0 \subset
W, \\
t \in \End W,\, t(D_i) \in D_{i-1}, \dim (D_i/D_{i-1}) = w_i, \im t
\subset V_0 \subset \ker t\}.
\end{split}
\end{equation*}
Consider the projection
\begin{equation*}
\begin{split}
\mathfrak{T}(w_1,\dots,w_r) \to \{\mathbf{D}=\{D_i\}_{i=0}^r \;|\;
0=D_0 \subset D_1 \subset \dots \subset D_r = W, \\
\dim(D_i/D_{i-1}) = w_i\}
\end{split}
\end{equation*}
given by $(\mathbf{D},V_0,t) \mapsto \mathbf{D}$.
It is easy to see that the fibers of this map are all
isomorphic and that in \cite{Savage} one could replace the tensor
product variety $\mathfrak{T}(w_1,\dots,w_r)$ by this fiber,
restrict the constructible functions to this fiber and the theory
would remain unchanged.  Let
$\mathfrak{T}_{\mathbf{D}}(w_1,\dots,w_r)$ denote the fiber over a
flag $\mathbf{D}$.  If we define
\[
D_i = W_0 \oplus \dots \oplus W_i,\; 0\le i \le r,
\]
then obviously
\[
\mathfrak{T}_{\mathbf{D}}(w_1,\dots,w_r) \cong \mathcal{L}(w_1,\ldots,w_r)
\]
and in the sequel we will identify these two varieties.

\paragraph{5.2.}
If $b \in CM_{w_1, \dots, w_r}^\mu$ is an unoriented
crossingless match, let 
\[
Y_b = \{(\mathbf{D},V_0,t) \in \mathfrak{T}(w_1, \dots, w_r)\; |\;
\dim (\ker t \cap D_i)/(\ker t \cap D_{i-1}) = b_i\}
\]
where $b_i$ is the number of left end points (of lower curves) and
lower end points (of middle curves) contained in the box representing $V_{w_i}$.
It is shown in \cite{Savage} (Proposition~3.2.1) that $\sqcup_b Y_b =
\mathfrak{T}(w_1, \dots, w_r)$ and that the closures of the $Y_b$ are
precisely the irreducible components of $\mathfrak{T}(w_1, \dots, w_r)$.
Let $X_b = Y_b \cap \mathcal{L}(w_1,\dots,w_r)$.
Then obviously
$\mathcal{L}(w_1, \dots, w_r) = \sqcup_{b \in LCM} X_b$.

\begin{prop}
$\mathcal{L}_l(w_1,\dots,w_r) = \sqcup_{b \in {^lLCM}} X_b$.
\end{prop}

\begin{proof}
We see from equation~\eqref{Nl_def} that $\mathcal{L}_l(w_1, \dots,
w_r)$ is the set of all $(t,V_0) \in \mathcal{L}(w_1, \dots, w_r)$
such that
\[
\dim \ker t|_{W_1 \oplus \dots \oplus W_i} \le l + \rank t|_{W_1
  \oplus \dots \oplus W_{i-1}} \; \forall \; 1 \le i \le r.
\]
Now, by the definition of the $X_b$, if $(t,V_0) \in X_b$ for some $b
\in LCM$ then $\dim \ker t|_{W_1 \oplus \dots \oplus W_i}$ is equal to
$\sum_{j=1}^i w_j$ minus the number of lower curves with both
end points among the lower $i$ boxes.  Also, $\rank t|_{W_1 \oplus
  \dots \oplus W_{i-1}}$ is equal to the number of lower curves with both
  end points among the lower $i-1$ boxes.  Let $c_i$ denote the number
  of curves with both end points among the lower $i$ boxes.  Then
\begin{gather*}
\dim \ker t|_{W_1 \oplus \dots \oplus W_i} \le l + \rank t|_{W_1
  \oplus \dots \oplus W_{i-1}} \\
\Leftrightarrow \sum_{j=1}^i w_j - c_i \le l + c_{i-1}
  \\
\Leftrightarrow \sum_{j=1}^i w_j - 2c_{i-1} - \#\{\text{curves with
  right end point in $i^{th}$ box}\} \le l
  \\
\Leftrightarrow 
\begin{split}
\sum_{i=1}^n w_i - \#\{\text{end points in first $i-1$ boxes of lower
  curves with both} \\
\text{end points in first $i$ boxes}\} \le l
\end{split}
\end{gather*}
and this is easily seen to be equivalent to the condition that $b \in {
  ^lLCM}$ (with the default bracketing).
\end{proof}

\paragraph{5.3.} We will now define a $U(\mathfrak{sl}_2)$-module structure on
a certain space of constructible functions on $\mathcal{L}_l(w_1, \ldots, 
w_r)$.
For $\mathbf{a} \in OLCM_{w_1, \dots, w_r}$, let ${\bar{\mathbf{a}}}$
be the associated element of $LCM_{w_1, \dots, w_r}$ obtained by forgetting the
orientation.  Define
\[
Y_{\mathbf{a}} = \{(\mathbf{D}, V_0, t) \in Y_{\bar{\mathbf{a}}}\; |\;
\dim W = \#\{\text{up-oriented vertices of $\mathbf{a}$}\}\}
\]
where the right end points of lower curves are oriented up (as well as
the up-oriented unmatched vertices).
Let $X_{\mathbf{a}} = Y_{\mathbf{a}} \cap
\mathcal{L}(w_1, \dots, w_r)$.  Then it follows from equation~(33) of
\cite{Savage} that
\[
X_b = \bigcup_{\mathbf{a} : \bar{\mathbf{a}} = b} X_{\mathbf a}.
\]
Now let
\[
\mathcal{B}^l_s = \{\id_{Y_{\mathbf{a}}}\; |\; \mathbf{a} \in {^lOLCM}\}
\]
where $\id_A$ is the function that is equal to one on the set $A$ and
zero elsewhere.  Let
\[
\mathcal{T}^l = \mathcal{T}^l_s(w_1, \dots, w_r) = \Span
\mathcal{B}^l_s.
\]
We endow $\mathcal{T}^l$ with the structure of a
$U(\mathfrak{sl_2})$-module as in \cite{Savage}.

\begin{theo}
$\mathcal{T}^l_s(w_1, \dots, w_r)$ is isomorphic as a
$U(\mathfrak{sl}_2)$-module to $V_{w_1} \otimes_l \dots \otimes_l
V_{w_r}$ and $\mathcal{B}^l_s$ is a basis for $\mathcal{T}^l_s(w_1, \dots, w_r)$
adapted to its decomposition into a direct sum of irreducible
representations.  That is, for a given $b \in {^lCM_{w_1, \dots,
    w_r}^\mu}$, the space $\Span \{\id_{Y_{\mathbf{a}}} \; |\;
\bar{\mathbf{a}} = b \}$ is isomorphic to the irreducible
  representation $V_\mu$ via the map
\[
\id_{Y_{\mathbf{a}}} \mapsto {^\mu v_{\mu - 2\#\{\text{unmatched down-oriented vertices of $\mathbf{a}$}\}}}.
\]
\end{theo}

\begin{proof}
The second part of the theorem follows from Theorem~3.3.1 of
\cite{Savage}.  Then
\begin{align*}
\mathcal{T}^l &= \bigoplus_{\mu} \bigoplus_{b \in {^lCM}_{w_1, \dots,
      w_r}^\mu} \Span \{\id_{Y_{\mathbf{a}}} \; |\; \bar{\mathbf{a}} =
      b\} \\
&\cong \bigoplus_{\mu} \bigoplus_{b \in {^lCM}_{w_1, \dots,
      w_r}^\mu} V_\mu \\
&\cong \bigoplus_{\mu} {^lH}_{w_1, \dots, w_r}^\mu \otimes V_\mu \\
&\cong V_{w_1} \otimes_l \dots \otimes_l V_{w_r}
\end{align*}
where ${^lH}_{w_1, \dots, w_r}^\mu$ is given the trivial module
structure. 
\end{proof}

\paragraph{Remarks.} We have used here the standard bracketing $(
\cdots (V_{w_1} \otimes_l V_{w_2}) \otimes_l V_{w_3}) \cdots \otimes_l
V_{w_r})$.  However, one could easily modify the definitions to use
any other bracketing.  The proofs would need only slight
changes.  Of course, as noted above, while we would still
recover the structure of the fusion product, the varieties involved
would be non-isomorphic in general.


\small{
}
\vspace{4mm}
\noindent
Olivier Schiffmann,\\
DMA ENS Ulm, 45 rue d'Ulm, 75005 Paris FRANCE;
email:\;\texttt{schiffma@dma.ens.fr}\\
Alistair Savage,\\
Department of Mathematics, Yale University, P.O. Box 208283, New
Haven, CT 06520-8283, USA;
email:\;\texttt{alistair.savage@aya.yale.edu}

\end{document}